\documentclass[a4paper]{amsart}
\usepackage{amsthm}
\usepackage{amssymb}
\usepackage{xy}
\usepackage{comment}
\xyoption{all}
\usepackage{float}
\usepackage{xcolor}
\usepackage{graphics}
\usepackage{amssymb}

\usepackage[utf8]{inputenc}
\usepackage[T1]{fontenc}
\usepackage{dsfont}
\usepackage[normalem]{ulem}
\usepackage{calligra}
\usepackage{verbatim}

\usepackage{color}

\definecolor{vermelho}{rgb}{0.9, 0.0, 0.0}

\newcommand{\Sss}{S}
\newcommand{\Ee}{E}
\newcommand{\ep}{\mathrm{\epsilon}}

 \newcommand{\EEE}{E_1,\dots,E_n}
 \newcommand{\lrarrow}{{\longrightarrow}}
  \newcommand{\Proj}{\mathrm{Proj}}
  \newcommand{\sgldim}{\mathrm{s.gl.dim}}
  \newcommand{\stgldim}{\mathrm{st.gl.dim}}
 
 \newcommand{\si}{\mathrm{\sigma}}

 \newcommand{\XX}{\mathbb{X}}
 \newcommand{\coh}[1]{\mathrm{coh}(#1)}

 \newcommand{\vc}{\vec{c}}

 \newcommand{\vom}{{\vec\omega}}
 
 \newcommand{\vx}{{\vec{x}}}
 \newcommand{\vy}{{\vec{y}}}
 
 \newcommand{\vz}{{\vec{z}}}

 \newcommand{\Hh}{\mathcal{H}}
 
 \renewcommand{\deg}{\mathrm{deg}\,}
 \newcommand{\ddim}{\mathrm{\dim}_k}
 \newcommand{\rk}{\mathrm{rk}\,}

 \newcommand{\Lp}{\mathbb{L}(\mathbf{p})}

 \newcommand{\Hom}{\mathrm{Hom}}
 \newcommand{\Ext}{\mathrm{Ext}}
 \newcommand{\HomX}{\mathrm{Hom}}

 \renewcommand{\phi}{\varphi}

 \newcommand{\Proof}{\par\noindent\textbf{Proof.}\ }
 \newcommand{\pp}{\mathbf{p}}
 \newcommand{\la}{\lambda}
 \newcommand{\ra}{\rightarrow}
 
 \newcommand{\ZZ}{\mathbb{Z}}
 
 \newcommand{\NN}{\mathbb{N}}
 \newcommand{\PP}{\mathbb{P}}
 \newcommand{\End}{\mathrm{End}}
 \newcommand{\Oo}{\mathcal{O}}
 \renewcommand{\mod}[1]{\mathrm{mod}(#1)}
 \def\ende{~\hfill$\square$}

 \newcommand{\Knull}{\mathrm{K}_0}
 
 \newcommand{\DHomX}{\mathrm{D\,Hom}_{\XX}}

 \newcommand{\ext}{{\rm Ext}^1}

 \newcommand{\ExtX}{\mathrm{Ext}^1}
 
 \newcommand{\DHom}{\mathrm{D\,Hom}}

 \newcommand{\lala}{{\mathbf \lambda}}

\definecolor{ashy}{gray}{0.6}

 \theoremstyle{plain}
 \newtheorem{Thm}{Theorem}[section]
 \newtheorem{Prop}[Thm]{Proposition}
 \newtheorem{Lem}[Thm]{Lemma}
 \newtheorem{Cor}[Thm]{Corollary}

 \theoremstyle{definition}
 \newtheorem{Defi}[Thm]{Definition}

 \newtheorem{Rem}[Thm]{Remark}

 \numberwithin{equation}{section}

\advance\textheight by 0.5cm
\advance\textwidth by 1.0cm
\advance\evensidemargin by -0.8cm
\advance\oddsidemargin by -0.8cm

\begin{document}

\title[Braid group action]{The braid group action on exceptional sequences for weighted projective lines}

\author[E. R. Alvares, E.N. Marcos and H. Meltzer]{Edson R. Alvares,
  Eduardo N. Marcos, and Hagen Meltzer}

\address{ Departamento de Matem\'atiica Universidade Federal do Parana, Brazil}

\email{rolo1rolo@gmail.com}

\address{ Departamento de Matem\'atica , IME , Universidade de Sao Paulo Brazil}

\email{enmarcos@ime.usp.br.com}

\address{Instytut Matematyki, Uniwersytet Szczeci\'nski, 70451 Szczecin, Poland}

\email{hagen.meltzer@usz.edu.pl}



\thanks{\footnotesize The third mentioned author thanks FAPESP,  from the grant 2019/08284-4, which made this work possible, the first mentioned author thanks FAPESP, from the grant 2018/08104-3.
The second  mentioned author was supported by the thematic project of FAPESP 2014/09310-5, a research grant from CNPq 302003/2018-5  and we all acknowledges support from the "Brazilian-French Network in Mathematics". }

\maketitle
\begin{center}
 We dedicate this work to the memory of Andrzeyj Skowro\'nski
\end{center}

\begin{abstract}
We give a new and intrinsic  proof of the transitivity of the braid group action
on the set of full exceptional sequences of coherent
sheaves on  a weighted projective line. We do not use here the corresponding result of Crawley-Boevey
 for modules over hereditary algebras.
As an application we prove that the strongest global dimension
of the category of coherent sheaves on a weighted projective line $\XX$
 does not depend on the parameters of $\XX$.
  Finally we prove that the determinant of the matrix obtained by taking the values of $n$ $\ZZ$-linear functions
  defined on the Grothendieck
group $K_0(\XX) \simeq \ZZ^n  $
of the elements of a full exceptional sequence is an invariant, up to sign.

\end{abstract}

\section{Introduction} \label{introduction}

Let $\XX$ be a weighted projective line in the sense of Geigle and Lenzing  \cite{Geigle:Lenzing:1987}.
The braid group $B_n$ on $n$ strings acts on the set
of full exceptional sequences in the category $\coh \XX$ of coherent sheaves on $\XX$, where $n$ denotes the rank
of the Grothendieck group $K_0(\XX)$ of $\coh \XX$. This action is given by mutations in
the sense of Gorodentsev and Rudakov \cite{GR}.
The following result was proved in \cite{Meltzer:1995}

\begin{Thm}\label{main}
The action of the braid group  on the set of full exceptional sequences in the category of
 coherent sheaves on a weighted projective line $\XX$ is transitive.
\end{Thm}

The proof was based on induction on the rank of the Grothendieck group of  $\coh \XX$ and on the rather strong result
of Crawley-Boevey \cite{CB} which states  that the braid group acts transitively on the set of full exceptional
 sequences in the category
of finitely generated modules over a hereditary algebra over an algebraically closed field.

It is desirable to have in the  geometric situation a purely sheaf-theoretical proof for the transitivity of the
 braid group operation. In this paper we show that this in fact can be done using perpendicular calculus
of exceptional pairs. For this we calculate the left perpendicular category of the sum of
 two line bundles  $L \oplus L(\vc)$ formed in the sheaf category, where $L$ is a line bundle and $\vc$ the canonical
 element of
the grading group of $\XX$.
For the convenience of the reader we also state the unchanged parts
of the original proof.

Furthermore, we give two applications of the transitivity of the braid group action.  First we show that the strongest
global dimension of a weighted projective
line  $\XX$,  a notion which we defined in this paper, is independent  of the parameters of $\XX$. This means that if $\XX=\XX(\pp,\lala)$ and
$\XX'=\XX(\pp,\lala')$,
are weighted projective lines with the same weight sequence $\pp$  and different parameter sequences
 $\lala$ and $\lala'$  then the strong global dimensions for $\XX$ and $\XX' $ are the same.

Second we prove that the determinant of the matrix obtained by applying $n$ additive functions
defined on the Grothendieck group of $\coh \XX$
  to the sheaves of a full exceptional sequence
on $\XX$ is independent of the exceptional sequence,  up to sign.
 Finally, we calculate this invariant for taking the rank function, the degree function and $n-2$ Euler forms with
 respect to simple exceptional sheaves.

\section{Preliminaries} \label{preliminaries}

\subsection {}
Weighted projective lines were introduced by Geigle and Lenzing in 1987  in
order to give a geometric approach to Ringel's canonical algebras \cite{Ri}. We recall some of the basic facts
and refer for details to \cite{Geigle:Lenzing:1987}.

Let $k$ be an algebraically closed field.
A weight sequence $\pp=(p_1, \dots, p_t)$ is a sequence of natural numbers $p_i$ with $p_i \geq 2$.
For a weight sequence $\pp$ denote by $\Lp$ the abelian group with generators $\vx_1, \dots, \vx_t$ and relations
$p_1 \vx_1= \dots = p_t \vx_t :=\vc$. The element $\vc$ is called the canonical element. $\Lp$ is an ordered group
with $\sum_{i=1}^t \NN \vx_i$ as cone of non-negative elements.
Furthermore, each element $\vx$ can be written, on a unique way,  in normal form  $\vx= l\vc + \sum_{i=1}^t l_i \vx_i $with
$l \in \ZZ$ and $0 \leq l_i < p_i$.
Consider further a sequence of parameters
$\lala=(\la_3, \dots, \la_t)$, that is the $\la_i $ are non-zero and pairwise distinct elements from $k$.
We denote $S=S(\pp,\lala)= K[X_1, \dots, X_t]/ (X_i^{p_i}+ X_1^{p_1} + \la_i  X_2^{p_2}, \, i=3,\dots ,t).$
The algebra $S(\pp,\lala)$ is $\Lp$-graded by defining $\deg (X_i)=\vx_i$.
Then the weighted projective line $\XX=\XX(\pp,\lala)$ is defined to be the $\Lp$-graded projective scheme
$\Proj^{\Lp} (S(\pp,\lala))$ and  the category $\coh \XX$ of coherent sheaves on $\XX$ is the quotient of
the category of finitely generated $\Lp$-graded $S$ modules modulo the $\Lp$-graded $S$  modules of finite length.
The category $\coh \XX$ is abelian, hereditary, that is $\Ext^i(A,B)=0$ for all $A$ and $B$ in
$\coh \XX$ and $i \geq 2$, and it has finite dimensional $\Hom$ and $\ext$ spaces. Moreover,
$\coh \XX$ admits Serre duality in the form $\ext(A,B) \simeq \DHom(B, A(\vom))$,
where $\vom$ denotes the dualizing element $(t-2)\vc -\sum_{i=1}^t  \vx_i$,
and consequently $\coh \XX$ has Auslander-Reiten sequences.

We denote the structure sheaf on $\XX$ by $\Oo$.
It is well known that the isomorphism class of line bundles on $\XX$  form a group, via the tensor product
and this group
is isomorphic to the group $\Lp$ via the
map $\vx \mapsto \Oo(\vx)$ where $\Oo(\vx)$ is the twisted by $\vx$ structure sheaf.
 Moreover, the homomorphism space between two line bundles can be calculated
by the formula $\Hom(\Oo(\vx),\Oo(\vy)) \simeq S_{\vy-\vx}$ and if $\vz= l \vc +\sum_{j=1}^t l_i\vx_i$ is
in normal form, then $\dim S_{\vz}=l+1$ provided $l \geq -1$.
For coherent sheaves on $\XX$ we have the rank and the degree function,
 which are defined on the Grothendieck group $K_0(\XX)$. The sheaves of rank $0$
are those of finite length.
One of the key results in \cite{Geigle:Lenzing:1987} is that the sheaf $\bigoplus_{0 \leq \vx \leq \vc} \Oo(\vx)$
 is a tilting sheaf
such that its endomorphism algebra is a canonical algebra.

\vspace{0,2cm}

\subsection{}\label{endsimple}
Recall that an object in a hereditary $k$-category $\Hh$ is called exceptional
 if $\End(E)=k$ and $\Ext^1(E,E)=0$. Moreover, a sequence of exceptional objects
$ \epsilon =(E_1, \dots, E_r)$ is called an exceptional sequence if $\Hom(E_j, E_i)=0=\Ext^1(E_j, E_i)$
for all $j>i$. If $r=2$ then $ \epsilon$ is called an exceptional pair and if $r$ equals the rank of the
Grothendieck group $K_0(\Hh)$ then $ \epsilon$ is called a full
exceptional sequence. This nomenclature is justified since every exceptional sequence has at most $K_0(\Hh)$  entries
and any exceptional sequence can be extended to at least one  full exceptional sequence.

Gorodentsev and Rudakov defined  mutations of exceptional sequences
 on $\PP^n$ which give rise
to an operation of the braid group
 $B_r= \langle \si_1, \dots , \si_{r-1} | \si_i \si_j =\si_j\si_i \;
 {\mathrm {for} } \; i-j\geq 2 \, \; {\mathrm {and} }\; \si_i \si_{i+1}\si_i=\si_{i+1}\si_i\si_{i+1} \rangle$
 on the set of (isomorphism classes) of exceptional sequences of length $r$ \cite{GR}.
For a categorical treatment we refer to \cite{B}.

We will study the action of the braid group $B_n$,
where $n$ is the rank of $K_0(\XX)$,
 on the set of full exceptional sequences in $\coh \XX$.
In this case each line bundle is exceptional. Moreover, the simple exceptional sheaves of rank $0$
fit in exact sequences

$$ 0 \lrarrow \Oo(j\vx_i) \lrarrow \Oo((j+1)\vx_i) \lrarrow S_{i,j} \lrarrow 0.$$

We continue this section  with the following lemma, which is probably well known, but for the sake of completeness we state and give a proof.

\begin{Lem}\label{elementary}
Let $\mathcal{A}$ be an abelian $k$-category and $M$ an object in it whose endomorphism ring is a finite dimension $k$-algebra, and $f$ an element in
$End({M})$. Then the following are equivalent:
\begin{enumerate}
    \item $f$ is a monomorphism,
    \item $f$ is an epimorphism
    \item $f$ is an isomorphism
\end{enumerate}
\end{Lem}

\Proof

We show that (1) implies (3).
Let $f$ be a monomorphism. We can assume that $f$ is non-zero. There is  $n$ such that $\{ f, f^2, \cdots, f^n \}$ is linearly dependent.  So there is a
non trivial linear combination
$\lambda_t f^t + \cdots + \lambda_n f^n = 0$ with $\lambda_t$ and  $\lambda_n$ non-zero.
So we have $f^t(\lambda_t Id + \cdots + \lambda _n f^{n -t}) = 0$. Since $f^t$ is a monomorphism, we get that $(\lambda_t Id + \lambda_{t+1} f + \cdots+
\lambda _n f^{n -t}) = 0$ which implies that
$\lambda_t Id =-  (\lambda_{t+ 1} f+\cdots +  \lambda _n f^{n -t})$.  Factoring out $f$  we get that $f$ is invertible, which shows that (1) implies (3.)
Analogously we have that (2) implies (3). Since clearly (3) implies (1) and (2), we have that the three assertions are equivalent.
\ende

\begin{Cor} {\label{corelementary}}  Let us assume  the same hypotheses as in lemma \ref{elementary}. If $M$ and $N$ are objects in $\mathcal{A}$ and there are
monomorphisms, $ f: M\to N$ and $g: N \to M$ then $f$ and $g$ are
isomorphisms. The analogous statement is valid for epimorphism. \ende
\end{Cor}

Given two sheaves  $A$ and $B$ over a weighted projective line, we define the trace map
${\mathrm {can}}: \Hom(A,B)\otimes A \to B $ in the usual way, i.e. ${\mathrm {can}}(f\otimes a) = f(a)$. In the literature, the image of ${\mathrm
{can}}$ is  also called the trace of $A$ in $B$.

Furthermore, if the space
$\Hom(A,B)$ is different from zero,  then the canonical map ${\mathrm {can}: \Hom(A,B) \otimes_k A \lrarrow B}$
 is surjective or injective but not bijective,
the proof for this fact is similar to the proof of \cite[Lemma
4.1]{HR}. In order to make our text complete,  we give it now.

\begin{Lem}
Let  $A, B$ be an exceptional pair in  $\coh \XX$, then the
trace map
${\mathrm {can}}: \Hom (A,B)\otimes A \to B$ is either a monomorphism or an epimorphism.
\end{Lem}
\Proof
We let $U$ be the image of ${\mathrm {can}}$,
 by $\mu$ the inclusion
$\mu: U \to B$ and ${\mathrm {can}} = \mu \delta$, where $\delta$ is induced by ${\mathrm {can}}$.

So we get the following exact sequence:

$$ (*)\ \ \ \ \  0\to U \to B \to B/U \to 0$$

Using that $\Ext^2 = 0$,
we get an epimorphism
$\Ext^1(B/U, \Hom(A,B)\otimes A) \rightarrow  \Ext^1(B/U, U)$.  This shows that the short
exact sequence $(*)$ comes from an extension in the group $\Ext^1(B/U, \Hom(A,B)\otimes A)$, i.e. we have the following  commutative diagram with exact
rows:

\begin{equation}{}
\xymatrix{
0    \ar[r]   &  \ar[d]^{\delta}  \Hom(A,B) \otimes A\  \ar[r]^{ \ \ \ \ \ \ \mu'}  & \ar[d]^{\delta'}  V  \ar[r] &  \ar[d]^{Id}  B/U \ar[r]  &            0\\
 0   \ar[r] &   U      \ar[r]^{\mu}  &  B \ar[r]    &  B/U \ar[r] &   0   }
\end{equation}

\noindent where $\delta$ and $\delta'$ are epimorphisms. From this diagram we get the following exact sequence:
$$0\to \Hom(A,B)\otimes A \stackrel{(\delta \ \mu')^{tr}}\longrightarrow U\oplus V \stackrel{(\mu \ -\delta')}\longrightarrow B\to 0.$$

Since $\Ext^1 (B, A) =0$ the exact sequence above splits.  We consider now two cases.

Case 1: $B$ is a direct summand of $U$. In this case,  there is  a monomorphism from $B$ to $U$, then since there is a monomorphism $\mu: U \rightarrow
B$, we use Corollary \ref{corelementary} and get that $\mu$ is an epimorphim. So the map ${\mathrm {can}}$ is an epimorphism.

Case 2:  $B$ is not a direct summand of $U$. Therefore, $U \simeq A^{t}$ for some $t$ and  $U \oplus V \simeq A^{t} \oplus (A^{s} \oplus B)$ where $t + s
= \dim \Hom(A,B)$. Since $U$ is the image of ${\mathrm {can}}$, we have  $\Hom(\Hom(A,B)\otimes A, B) = \Hom(\Hom(A,B)\otimes A, U) \simeq \Hom
(\Hom(A,B)\otimes A, A^{t})$ and $\dim \Hom(A,B) = t \dim \Hom(A,A)$. Therefore, $t = \dim \Hom(A,B)$ and $s = 0$.

The morphism  $\Hom(A,B)\otimes A \stackrel{\delta}\rightarrow  U$,  the isomorphism $U \simeq \Hom(A,B)\otimes A$ and the Corollary \ref{corelementary}
give us that $\delta$ is a monomorphism.
\ende

\vspace{0,2cm}

The left mutation of $(A,B)$ is the exceptional pair $(L_A B, A)$, where $ L_A B$ is given by one of
the following three exact sequences: if $ \Hom(A,B)\neq 0$ then

$$   0 \lrarrow L_A B  \lrarrow  \Hom(A,B) \otimes_k A \stackrel{{\mathrm {can}}}{\lrarrow} B \lrarrow 0,$$

$$   0 \lrarrow  \Hom(A,B) \otimes_k A \stackrel{{\mathrm {can}}}{\lrarrow} B  \lrarrow L_A B  \lrarrow 0,$$

\noindent
and if $ \ext(A,B)\neq 0$
then
$$   0 \lrarrow  B   \lrarrow L_A B   \lrarrow  \ext(A,B) \otimes_k  A \lrarrow 0, $$
where  the third sequence is the universal extension. If $\Hom(A,B)=0=\ext(A,B)$ then $L_A B=B$ and
the left mutation of the pair $(A,B)$ is called a transposition.
Now, the generators of $ B_{n}$ act on the set of full exceptional sequences in $\coh \XX$ as follows:
$$\si_i \cdot (E_1, \dots E_{i-1}, E_{i},  E_{i+1},  E_{i+2} ,\dots,  E_{n})=
(E_1, \dots E_{i-1}, L_{E_i} E_{i+1},E_i,  E_{i+2} ,\dots,  E_{n}).$$

Further the right mutation of an exceptional pair $(A,B)$ is the exceptional pair $(B, R_B A)$, where
 $ R_B A$ is given by one of
the following three exact sequences

$$   0 \lrarrow A   \stackrel{{\mathrm {cocan}}}{\lrarrow} { \mathrm {DHom}}(A,B) \otimes_k B \lrarrow  R_B A\lrarrow
0,$$

$$   0 \lrarrow R_B A  \lrarrow  A   \stackrel{{\mathrm {cocan}}}{\lrarrow} { \mathrm {DHom}}(A,B) \otimes_k B \lrarrow
0,$$

$$   0 \lrarrow   { \mathrm{ DExt}}^1(A,B) \otimes_k  B     \lrarrow  R_B A  \lrarrow    A \lrarrow 0, $$
where $ { \mathrm D}=\Hom_k(-,k)$, ${\mathrm {cocan }}$ denotes the co-canonical map and  the third sequence is
the universal extension.
Then $\si_i^{-1}$ acts in the following way.

 $$\si_i^{-1} \cdot (E_1, \dots E_{i-1}, E_{i},  E_{i+1},  E_{i+2} ,\dots,  E_{n})=
(E_1, \dots E_{i-1},E_{i+1},  R_{E_i+1} E_{i},  E_{i+2} ,\dots,  E_{n}).$$

The following lemma is a useful tool.

 \medskip
 \noindent

\begin{Lem}\label{helix}
We have

(i) $\si_1 \dots \si_{n-1} (E_1, E_2, \dots, E_n) = ( E_n(\vom),E_1, E_2, \dots, E_{n-1})$

(ii) $\si_{n-1} \dots \si_{1} (E_1, E_2, \dots, E_n) = (E_2, \dots, E_{n-1}, E_1(-\vom))$

(iii) In the orbit of an exceptional sequence $(E_1,\dots E_a, E_{a+1}, \dots )$ there
is an exceptional sequence of the form $ ( E_a, E_{a+1},\dots)$
\end{Lem}
\medskip
The proof for (i) and (ii) is given in \cite[Proposition 2.4]{Meltzer:1995} and
(iii) is a consequence of (i) and (ii).

\medskip
\subsection{}
Recall that for an object $X$ in a hereditary category  $\Hh$ the left perpendicular category
with respect to $X$ is defined as the full subcategory of all objects $Y$ satisfying
$\Hom(Y,X)=0$ and $\ext(Y,X)=0$ (see \cite{Geigle:Lenzing:1991}). The right perpendicular category is defined dually.

\section{Proof of Theorem 1.1 } \label{proof}
\subsection{}

In this section we will prove Theorem \ref{main}.
Let $\XX$ be a weighted projective line of weight type $\pp=(p_1, \dots, p_t)$ and rank of
$K_0(\XX)$  equals $n$.
 We start with the following observation.

\begin{Prop}\label{perpendicular}

(a) Let $(L,L')$ be an exceptional pair of line bundles in  $\coh \XX$ with $\dim_k\Hom(L,L')\geq 2$.
Then $L' \simeq L(\vc)$ and  $\dim\Hom(L,L')=2$.

(b) The left perpendicular category with respect to $L \oplus  L(\vc)$ for a line bundle $L$, formed in $ \coh \XX $,
 consists only of finite length sheaves.
Moreover,
this perpendicular category is
 equivalent to the category of finite dimensional modules over the path
 algebra of the disjoint union  of linearly oriented quivers of type $A_{p_i-1}$, $i=1, \dots, t$.
\end{Prop}

\Proof
(a)
We have $L'=L(\vx)$ for some $\vx$. We write  $\vx $  in  normal form and, after renumbering
the indices, if necessary,
 $\vx =l \vc +\sum_{j=1}^r l_j\vx_j$,
where $l_{1}\neq 0, \dots, l_{r} \neq 0$ for some $r$. Since $\dim\Hom_k(L,L') \geq 2$
we have $l \geq 1$. Using Serre duality and
the fact that $(L,L(\vx))$ is an exceptional pair we have
$0=\Ext^1(L(\vx), L )\simeq \Hom(L,L(\vx+\vom))\simeq \Hom(\Oo,\Oo(\vx+\vom))   $. Now
$\vx+\vom=  l \vc +\sum_{j=1}^r l_j \vx_j +   (t-2)\vc -\sum_{i=1}^t \vx_i
                   =  (l-2+r) \vc +\sum_{j=1}^r (l_j-1)\vx_j  +\sum_{i=r+1}^t (p_i-1) \vx_i$.
This element is in normal form and it follows that $l-2+r <0$, hence $l=1$ and $r=0$.
Consequently   $\vx=\vc $.

(b) After  renumbering
 the indices, if necessary,  for the simple exceptional sheaves in the tubes we can assume that
$ \Ext^1(S_{i,0}, L ) \neq 0$ for $i=1,\dots, t$.

The Riemann-Roch formula \cite[2.9]{Geigle:Lenzing:1987}  applied to $S_{i,0}$ and $L$ yields

 $$\sum_{j=0}^{p-1}\langle  \tau^j S_{i,0}, L\rangle=
 p\, (1-g)\, \rk (S_{i,0}) \, \rk( L) +
\det \begin{pmatrix} \rk  (S_{i,0})&  \rk (L )\\ \deg (S_{i,0} )&  \deg ( L )\end{pmatrix}$$

\noindent
where $p$ denotes the least common multiple of the weights $p_1,\dots,p_t$,
 $g$ is the genus of the weighted projective line and
 $\langle A,B \rangle= \dim \Hom(A,B)- \dim \ext(A,B)$ the Euler form.
Since the $\tau$ period of $S_{i,0}$ is $p_i$, $\rk(S_{i,0}) =0$ and
 $\deg(S_{i,0}) =\frac{p}{p_i}$ we conclude that
$\frac{p}{p_i} \sum_{j=0}^{p_i-1} \langle  \tau^j S_{i,0}, L
  \rangle =-\frac{p}{p_i}$.

Since there
are no non-zero homomorphisms from finite length sheaves to vector bundles we obtain that
$
  \sum_{j=0}^{p_i-1}\dim\ext ( \tau^j S_{i,0},L)=1$
 and therefore
$\ext (\tau^j S_{i,0},L)=0$ for  $ j=1,\dots, p_i-1$. Again using that there
are no non-zero homomorphisms from finite length sheaves to vector bundles we see
that the sheaves
$S_{i,j}$ for $ j=1,\dots, p_i-1$, $i=1,\dots, t$ belong to the left perpendicular category, formed in $\coh \XX$,
$^{\perp}(L)$,formed in $\coh \XX$.
 Since $S_{i,j}=  S_{i,j}(\vc)$, the same argument can be applied to the
line bundle $L(\vc)$.
Therefore the finite length sheaves $S_{i,j}$ for $ j=1,\dots, p_i-1$, $i=1,\dots, t$ belong
to the left perpendicular category
 $ \Hh = ^{\perp}\!(L \oplus L(\vc))$.

 The category $\Hh$ can be obtained by forming first the left perpendicular category $\Hh_1$
 with respect to $L$ in $\coh \XX$ and then the left   perpendicular
 category
 $\Hh_2$
 with respect to $L(\vc)$ in $\Hh_1$.
 The category $\Hh_1$ is known to be equivalent to the category of finitely generated modules over
  a hereditary algebra, in fact the path algebra of the quiver, obtained from the
quiver of the canonical algebra $\End( \bigoplus_{0 \leq \vx \leq \vc} L(\vx))$
 by removing the  vertex which corresponds to
$L$, see ( \cite{LP}).
Then by a result of Happel \cite{Ha},
using the fact that $L(\vc)$ considered in the module category $\mod {H_1}$ is exceptional,
the category $\Hh_2$ is  equivalent to the category of finitely generated modules over
a hereditary algebra $H_2$ . Moreover, both results  together imply that the rank of the Grothendieck group
 $K_0(\Hh_2)$ equals $n-2$.

 Denote by  $^{[j]} S_{i,1}$ the indecomposable sheaf with socle $ S_{i,1}$ and
quasi-length $j$.
The sheaf $ T = \bigoplus_{i=1}^{t}  \bigoplus_{j=1}^{p_j-1} \, ^{[j]}S_{i,1}  $
satisfies $\Ext^1(T,T)=0$ and consists of $n-2$ indecomposable direct summands.
Therefore $T$
is a tilting sheaf in  $\Hh$ and consequently $\Hh$  consists of the objects of the wings
for $^{[j]} S_{i,p_i-1}$, $i=1, \dots t$.
This  shows that the endomorphism algebra of $T$ is the disjoint union  of linear quivers
of type $A_{p_i-1}$, $i=1, \dots, t$.
\ende

\medskip

\subsection{}
We will use the following three results of \cite{Meltzer:1995}

\begin{Lem}    \cite[Lemma 2.7]{Meltzer:1995}
        \label{unique}
       Two distinct complete  exceptional sequences, differs in at least two
       places.
\end{Lem}

\medskip

{\bf 3.5}
An exceptional sequence $(E_1,...,E_n)$ in $\coh \XX$ is called orthogonal if
$\Hom(E_i,E_j)=0$ for all $i \neq j$.

\begin{Prop}\cite[Proposition 2.8]{Meltzer:1995} \label{orthogonal}
There are no orthogonal complete exceptional sequences in $\coh \XX$.
\end{Prop}

\begin{Lem}\cite[Lemma 3.1]{Meltzer:1995}\label{monoepi}
Let $\EEE$ be an exceptional sequence in $\coh \XX$ such that
$\dim_k\HomX(E_1,E_2) \ge 2$.

(i)
 Suppose that $LE_2=L_{E_1}E_2$ is defined by an exact sequence
$$
0 \ra  LE_2 \ra \HomX(E_1,E_2)\otimes E_1 \ra E_2 \ra 0.
$$
\noindent Then morphisms $ 0 \neq h \in \HomX(LE_2,E_1)$ and $0 \neq f \in \HomX(E_1,E_2)$
 are either  both monomorphisms or both epimorphisms.

(ii)
 Suppose that $RE_1=R_{E_2}E_1$ is defined by an exact sequence
$$
0 \ra E_1 \ra \DHomX(E_1,E_2) \otimes E_2 \ra RE_1 \ra 0.
$$
\noindent Then morphisms $ 0 \neq h \in \HomX(E_2,RE_1)$ and $0 \neq f \in \HomX(E_1,E_2)$
 are either both monomorphisms or both epimorphisms.
\end{Lem}

\medskip
\subsection{}
For an  exceptional  sequence ${\ep}= (\EEE)$  we define
$$\|\ep\|=\\(\rk (E_{\pi(1)}),...,\rk (E_{\pi(n)})),$$
 where $\pi$ is a permutation of $1,...,n$
such that $\rk (E_{\pi(1)}) \ge ...\ge \rk (E_{\pi(n)})$ .

\begin{Prop}\label{reduction}
Let $\XX$ be a weighted projective line with at least one weight, i.e. $\XX \neq
 \PP^1$. Then in  each
orbit, under the braid group action, there is a complete exceptional sequence containing a simple sheaf of
rank $0$.
\end{Prop}

\Proof
We show first the following claim:
 if $\epsilon=E_1, \dots, E_n)$ is a  complete exceptional sequence in $\coh\XX$ with $\rk (E_i)\geq 1$
for all $i$
then there exists $\si \in B_n$
such that $\|\si \cdot \epsilon\| < \| \epsilon \|$.

Let $\epsilon=(E_1, \dots, E_n)$
 be a  complete exceptional sequence in $\coh\XX$
with $\rk (E_i)\geq 1$
for all $i$.
We know from \ref{orthogonal} that $\ep$ is not orthogonal.
Choose $a<b$ such that $\Hom(\Ee_a,\Ee_b)\neq 0$, but  $\Hom(\Ee_i,\Ee_j) = 0$
for the remaining $a\le i < j \le b$.

Let $f: \Ee_a \ra \Ee_b$ a nonzero morphism.
 We know that  $f$ is a monomorphism or an epimorphism, thus we distinguish two cases.

{Case 1}: $f$ is a monomorphism.

\noindent
Then $f$ induces epimorphisms
$\ExtX(\Ee_b,\Ee_i) \twoheadrightarrow \ExtX(\Ee_a,\Ee_i)$
for all $i$. Since the first $\Ext$-group is zero for $i\le b$  the second
$\Ext$-group also vanishes for these $i$.
We see that both $\HomX(\Ee_a,\Ee_i)=0$, and $\ExtX(\Ee_a,\Ee_i)=0$ for all $a<i<b$,
therefore applying transpositions we obtain that
$$\si_{b-2}^{-1}...\si_{a+1}^{-1}\si_a^{-1} \ep = (\Ee_1,...,\Ee_{a-1},\Ee_{a+1},...,\Ee_{b-1},\Ee_a,
\Ee_b,...,\Ee_n).$$
Moreover, using  Lemma   \ref{helix},
we can assume that  $a=1$ and $b=2$.

Now, the left mutation $L\Ee_2=L_{\Ee_1}\Ee_2$ is defined by an exact sequence being
of the form
$$ {\rm (i)} \hspace{2cm} 0 \ra \HomX(\Ee_1,\Ee_2)\otimes\Ee_1 \ra \Ee_2 \ra  L\Ee_2 \ra 0 $$
or
$$ {\rm (ii)} \hspace{2cm}  0  \ra  L\Ee_2 \ra  \HomX(\Ee_1,\Ee_2)\otimes\Ee_1 \ra \Ee_2 \ra 0.$$

In the case (i) we have  $\rk (L\Ee_2) < \rk( \Ee_2) $, hence
$\| \si_1\ep\| < \|\ep\|$ and we are done.

In the case (ii) there exists a nonzero morphism $h: L\Ee_2 \ra \Ee_1$. Again, $h$ is a
monomorphism or an epimorphism. Because $f$ is a monomorphism we infer
 from the sequence (ii) that
$\dim_k\HomX(E_1,E_2)\ge 2$. But then, in view of Lemma \ref{Le3.1},
$h$ is a monomorphism.
 Thus
$$\rk (L\Ee_2) \le \rk (\Ee_1 )\le \rk(\Ee_2).$$
If $\rk (L\Ee_2 )< \rk (\Ee_2)$  we  apply $\si_1^{-1}$ as before
and obtain $\| \si_1^{-1}\ep\| < \|\ep\|$.

Assume otherwise that $\rk (L\Ee_2) = \rk (\Ee_2)$. Then also  $\rk (\Ee_1) = \rk (\Ee_2)$
and therefore $\dim_k\HomX(\Ee_1,\Ee_2)=2$.

Consider an exact sequence
$$
0 \ra E_1 \stackrel{f}{\ra} E_2 \ra C \ra 0$$
where $C={\rm coker (f)}$. Clearly $\rk (C)=0$. Furthermore,
applying the functor $\HomX(E_i,-)$  we conclude that $\dim_k\HomX(E_i,C)=1$, for $i=1,2$.
Finally, applying the  functor $\HomX(-,C)$
we obtain $\Hom(C,C)=k$ and $\ext(C,C)=k$, in particular $C$ is indecomposable.

 We have to consider two cases. First assume that
 C is a finite length sheaf  concentrated at an ordinary point.
 Now $\End(C)=k$ which implies that $C$ is a simple sheaf.
 The Riemann-Roch formula yields $\HomX(L,C)=k$ for each line bundle $L$.
  Thus using a line bundle filtration for
  $E_1$ we obtain $\dim_k\HomX(E_1,C)=\rk ( E_1)$. We have shown  before that
$\dim_k\HomX(E_1,C)=1$. Thus we obtain $\rk (E_1)=\rk (E_2)=1$ and we
have also $\dim\Hom(E_1, E_2)=2$.
 But then we have $\rk (E_i)=0$ for $i >2$ by Proposition
 \ref{perpendicular}.

Now, assume that $C$ is a sheaf of finite length concentrated at an
exceptional point, say  $\la_i$ of weight $p_i$.
It follows from $\Hom(C,C)=k$, $\ext(C,C)=k$  and the tube structure of the
Auslander-Reiten quiver  that the length of $C$ is $p_i$, and therefore
for the classes in the Grothendieck group $K_0(\XX)$ group we have
$[C]=\sum_{j=0}^{p_i-1} [\Sss_{i,j}]$ where $\Sss_{i,j}$ are the objects on the mouth of the tube.
From the exact sequences stated in subsection \ref{endsimple} we infer that
$ [\Sss_{i,j}] = [\Oo(j+1)\vx_i) ]  -[\Oo(j)\vx_i] $ for
 $i=1, \dots, t$, $j=1, \dots p_i$.
Hence $  \sum_{j=0}^{p_i-1} [\Sss_{i,j}] = [\Oo(\vc)]- [\Oo]$.
On the other hand there is an exact sequence
$0 \lrarrow \Oo  \lrarrow \Oo(\vc) \lrarrow \Sss \lrarrow 0 $
where $\Sss$ is a simple finite length sheaf concentrated in an ordinary point and consequently
$[C]=[\Sss]$.
 We conclude that
$1=\dim_k\HomX(E_1,C)= \chi ([E_1],[C])=\chi([E_1],[\Sss])=\dim_k\HomX(E_1,\Sss)=\rk (E_1)$,
where $\chi$ is the Euler form. Then we have $\rk (E_1)=\rk (E_2) =1$ and
$\dim\Hom(E_1,E_2)=2$, and consequently
$\rk (E_i)=0$ for $i >2$ by Lemma \ref{perpendicular},  contrary to our assumption.

{Case 2:} $f$ is an epimorphism.

Then $f$ induces  epimorphisms
$ \ExtX(E_i,E_a) \twoheadrightarrow \ExtX(E_i,E_b)$, for all $i$.
The first $\Ext$-group  is zero for $i \ge a$, thus also the second $\Ext$-group vanishes
 for  these $i$. We see that both $\HomX(E_i,E_b)=0$ and $\ExtX(E_i,E_b)=0$ for  all $a<i<b$,
 and again applying transpositions we have
$$\si_{a+1}^{-1}...\si_{b-1}^{-1}\ep=(E_1,...,E_{a-1},E_a,E_b,E_{a+1},...,E_n). $$
As before we can assume $a=1$ and $b=2$.
Then $RE_1=R_{E_2}E_1$ is defined by an exact sequence
$$
{\rm (i)} \hspace{2cm} 0 \ra RE_1 \ra E_1 \ra \DHomX(E_1,E_2) \otimes E_2 \ra 0 $$
or
$$ {\rm (ii) }\hspace{2cm} 0 \ra E_1 \ra  \DHomX(E_1,E_2) \otimes E_2 \ra RE_1 \ra 0
$$

In the first case  we have $\rk (RE_1) < \rk (E_1)$, and consequently
$\|  \si_1^{-1}  \ep\| < \|\ep\|$.
In the second case  there is a  nonzero map $h: E_2 \ra RE_1$, which again
is a monomorphism or an epimorphism.
Since $f$ is an epimorphism we conclude that $\Hom(E_1,E_2)\geq 2 $ and
therefore $h$ is an epimorphism by Lemma \ref{monoepi}.

Now, in this case,
$$ \rk( E_1) > \rk (E_2) > \rk (RE_1)$$
\noindent and therefore again
$\| \si^{-1}\ep\| < \|\ep\|$.

So the claim is proved.
We see that after applying  successively the norm reduction above, if necessary,
 we can shift by a braid group element
 any full exceptional sequence
to a sequence containing an exceptional sheaf of rank $0$.

 We will show now that
in the same orbit there is an exceptional sequence containing a  simple sheaf.

Now let $s$ be  the minimal number with the property that the orbit of ${\ep}$ contains an
exceptional sequence with  a  rank $0$ sheaf $F$ of length $s$.
By  Lemma \ref{helix}  we can assume that this exceptional sequence is of the form
$(\Ee_1,...,\Ee_{n-1},F)$.

We have to show that $s=1$. Assume contrary that $F$ is not simple and denote by
$\Sss$ the socle of ${F}$.
We claim that
$(\Ee_1,...,\Ee_{n-1},\Sss)$ is an exceptional sequence, too.
Indeed, we have $\Ext^{1}(\Sss,\Ee_i)=0$ for $1 \le i \le n-1$,
because the embedding $\Sss \hookrightarrow F$ induces epimorphisms
$\Ext^{1}(F,\Ee_i) \twoheadrightarrow \Ext^{1}(\Sss,\Ee_i)$ and the first $\Ext$-group vanishes  by assumption.
On the other hand,  $\Hom(\Sss,\Ee_i)=0$ for $1 \le i \le n-1$,  because the existence of a
nonzero morphism from $\Sss$ to some $\Ee_i$ implies that $E_i$ also has
 finite length,
  and equals therefore  $^{[r]}\Sss $,  for some $r$,  the unique indecomposable finite length
sheaf with socle $\Sss$ and length  $r$. Then $r \ge s$ by minimality of $s$.
But this  implies  $\Hom(F,\Ee_i)\neq 0$, contrary to the fact that $(\Ee_1,...,\Ee_{n-1},F)$
is an exceptional sequence.
Thus we have two exceptional sequences which coincide  in the first $n-1$ terms but
are different in the last one.
By Lemma \ref{unique}
 this is impossible.\ende

 \medskip

\subsection  \, Proof of Theorem\ref{main}
We show by induction on the rank $n$ of $\Knull(\XX)$ that the group $B_n$ acts transitively on
the set of complete exceptional  sequences in $\coh \XX$.

If $n=2$ then $\XX=\PP^1$. In this case  an exceptional  sequence is of the form $(\Oo(i),\Oo(i+1))$
for some $i \in \ZZ$
and the braid group $B_2 \cong \ZZ$ obviously acts transitively on the set of
these exceptional  sequences.

Now, suppose $n>2$ and assume that ${\ep}=(\EEE)$ is a full exceptional sequence in
$\coh \XX$. By Proposition  ~\ref{reduction} and applying, if necessary,
Lemma   \ref{helix}
 we have  $g.{\ep}=(E'_1,...E'_{n-1},\Sss)$ for some
$g \in B_n$ and some simple exceptional finite length sheaf $\Sss$. Denote by
${\kappa}=(\Oo,\Oo(\vx_1),\Oo(2\vx_1),\dots, \Oo((p_1-1)\vx_1) , \Oo(\vx_2), \dots \Oo((p_t-1)\vx_t),\Oo(\vc))$
the exceptional sequence corresponding to the canonical tilting sheaf.
Since $\Sss$ is exceptional simple,
$\Sss= \Sss_{i,j}  $ for some $i,j$.
From the exact sequence
$$
0 \ra \Oo(j\vx_i) \ra \Oo((j+1)\vx_i) \ra \Sss_{i,j} \ra 0
$$
we  see that the right mutation of the pair $(\Oo(j\vx_i),\Oo((j+1)\vx_i)$
equals $(\Oo((j+1)\vx_i),S_{i,j})$.
Thus,  for some $g_1 \in B_n$ we get
$g_1.{ \kappa}=(\Oo,...,\Oo((j+1)\vx_i),S_{i,j},...)$.
Observe that in case $j=p_i-1$, $i\neq t$, we first can apply transpositions
in order to arrange that $\Oo(j\vx_i)$ and $\Oo((j+1)\vx_i)$ are neighbours.
Applying if necessary,  Lemma 2.2, we obtain
$g_2.{ \kappa}=(F_1,...,F_{n-1},\Sss)$ for some $g_2\in B_n$ and line bundles $F_1,...,F_{n-1}$.
Now, by \cite{Geigle:Lenzing:1991} the right perpendicular category $\Sss^{\perp}$ is equivalent to a sheaf category
 $\coh {\XX'}$ for  a weighted   projective  line $\XX'=\XX({\bf p'},\lala)$ with
 weight sequence
${\bf p'}=(p_1,...p_{i-1},p_i-1,p_{i+1},...,p_n)$. By induction
$(E'_1,...,E'_{n-1})$ and $(F_1,...,F_{n-1})$, considered as complete exceptional sequences in
 $ \Sss ^{\perp}$, are in the same orbit under the action of the braid group  $B_{n-1}$
on the set of complete exceptional   sequences in $\coh {\XX'}$.
 We conclude that ${\ep}$ and ${\kappa}$ are in the same orbit,
which finishes  the proof.\ende

\section{The  strongest global dimension of $\coh \XX$} \label{strong}

\subsection{}

As an application of the transitivity of the braid group action in
\cite{Meltzer:2004} it was shown that several data are independent
of the parameters of the weighted projective line.
Here we are going to  investigate  the strongest global dimension of  weighted projective lines
by    studying  the spreading of tilting complexes in the derived category.

\begin{Defi}
   The strongest global dimension is the maximum of the strong global dimension of all algebras which
    are derived equivalent do $\coh \XX$
\end{Defi}

We have the following characterization of the strongest global dimension.

\begin{Prop}
The strongest global dimension of a weighted projective line
$\XX$ is  one if $\XX=\PP^1$ or
is
the maximal number $m+2$ such that there exists a
tilting complex $T$ in the derived category of $\coh \XX$ of the form
 $\bigoplus_{i=0}^m T_i [i]$ with $T_i \in \coh\XX$, $i \in \ZZ$ and $T_0 \neq 0 \neq T_m$.
\end{Prop}
\Proof See Theorem 1 in \cite{ALM}. \ende
\vspace{0.2 cm}

The  strongest global dimension of $\XX$ will be denoted  by
 $\stgldim \XX$.
 It follows from the definition that if the bounded  derived category of an algebra $A$ is
  triangular equivalent
   to the bounded derived category  of  $\coh \XX$, then
  the $\sgldim A \leq \stgldim \XX$.

In \cite{Meltzer:2004} it was shown that if $\XX$ has weight type $(2,2,\dots2)$, ($t$ entries)
 then  $\stgldim \XX=2$.
For tubular weighted projective lines $\XX$ bounds  of the strongest global dimension were
 given in \cite{S}.

 Before the main theorem in this section we have the following remarks and facts.
  Recall that, if $(A,B)$ is an exceptional pair in  $\coh \XX$, then $\Hom_{\XX}(A,B) = 0$ or $  \Ext^1(A,B) = 0$.
First we  assume that $\Hom(A, B)\neq 0$. We have then  two cases for the left mutation to consider:

$$(\alpha): 0 \lrarrow L_A B  \lrarrow  \Hom(A,B) \otimes_k A \stackrel{{\mathrm {can}}}{\lrarrow} B \lrarrow 0,$$
$$(\beta):  0 \lrarrow  \Hom(A,B) \otimes_k A \stackrel{{\mathrm {can}}}{\lrarrow} B  \lrarrow L_A B  \lrarrow 0.$$

It is important to note that the surjectivity of the canonical  map  depends only on $\rk (A)$, $\rk( B)$ and on the
dimension of the spaces $\Hom_{\XX}(A,B)$. We have that
\begin{center}
\begin{itemize}
    \item  if $\rk  (A) \neq 0$ then\\
    $\mathrm{can}$ is surjective $\Longleftrightarrow$ $\ddim  \Hom(A, B) \cdot rank \ A > rank \ B$.
    \item if $\rk  (A) = 0$ \\
    $\mathrm{can}$ is surjective $\Longleftrightarrow$ $\ddim  \Hom(A, B) \cdot \ddim \ A > \ddim  \ B$.
\end{itemize}

\end{center}

Applying $\Hom_{\XX}(\ , A)$ in $(\alpha)$ we have $\Hom_{\XX}(\alpha , A):$

\[
  \begin{array}{rcll}
    0 \lrarrow  \Hom(B  , A) &  \lrarrow   \Hom(\Hom(A,B) \otimes_k A , A)  & \lrarrow   \Hom(L_A B  ,
    A)  \lrarrow  \\
    \lrarrow   \Ext^1(B  , A)  & \lrarrow
  \Ext^1(\Hom(A,B) \otimes_k A , A) & \lrarrow     \Ext^1(L_A B  , A)   \lrarrow 0 \\
  \end{array}
  \]

Applying $\Hom_{\XX}(\ , A)$ in $(\beta)$ we have $\Hom_{\XX}(\beta , A):$
\[
\begin{array}{rcll}
    0 \lrarrow   \Hom(L_A B  , A) &  \lrarrow  \Hom(B  , A)  &  \lrarrow   \Hom(\Hom(A,B) \otimes_k A ,
    A) \lrarrow  \\
    \lrarrow   \Ext^1(L_A B  , A)  & \lrarrow
  \Ext^1(B  , A)  & \lrarrow     \Ext^1(\Hom(A,B) \otimes_k A , A)   \lrarrow  \\
  \end{array}
  \]

Now, the following remarks follows from both long exact sequences:
\begin{Rem} \label{dimension} The mutation of $(A, B)$ is the exceptional pair $(L_A B, A)$ and:
\begin{itemize}
 \item In the case  $(\alpha)$, the conditions  $\Hom_{\XX}(B  , A) = 0 =   \Ext^1(B  , A)$
  imply that \\
  $\ddim  \Hom(L_A B  , A) = \ddim  \Hom(A, B)$ and
 $\ddim   \Ext^1(L_A B  , A) = 0.$
 \item  In the case  $(\beta)$, the conditions  $\Hom_{\XX}(B  , A) = 0 =   \Ext^1(B  , A)$
  imply that \\
  $\ddim  \Hom(L_A B  , A) = 0$ and
 $\ddim   \Ext^1(L_A B  , A) = \ddim  \Hom(A, B).$
  \end{itemize}

Now assume that  $ \Ext(A,B)\neq 0$. Then have the universal extension
$$(\gamma):  0 \lrarrow  B   \lrarrow L_A B   \lrarrow  \Ext(A,B) \otimes_k  A \lrarrow 0. $$

Applying $\Hom_{\XX}(\ , A)$ in $(\gamma)$ we have $\Hom_{\XX}(\gamma , A):$
\[
\begin{array}{rcll}
    0 \lrarrow   \Hom(\Ext(A,B) \otimes_k  A, A) &  \lrarrow  \Hom(L_A B , A)  &  \lrarrow   \Hom(B, A)
    \lrarrow  \\
    \lrarrow   \Ext^1(\Ext(A,B) \otimes_k  A, A)  & \lrarrow
  \Ext^1(L_A B , A)  & \lrarrow     \Ext^1(B, A)   \lrarrow & 0. \\
  \end{array}
  \]
\end{Rem}

 \label{dimensionext}  The mutation of $(A, B)$ is the exceptional pair $(L_A B, A)$ and:
\begin{itemize}
 \item In the case  $(\gamma)$, the conditions  $\Hom_{\XX}(B  , A) = 0 =   \Ext^1(B  , A)$
  imply that \\
  $\ddim  \Hom(L_A B  , A) = \ddim Ext^{1}_{\XX}(A, B)$ and
 $\ddim   \Ext^1(L_A B  , A) = 0.$
 \end{itemize}

 \subsection{}

Summarizing, from \ref{dimension} and \ref{dimensionext}, if $(A. B)$ is an exceptional pair,  then on the  mutation
pair  $(L_A B, A)$  we can compute the dimensions
 $\ddim  \Hom((L_A B, A) $, $ \ddim  \Ext^1(L_A B, A) $, and $\rk (L_A B)$ without using the parameters
 $\lambda$.

\begin{Lem}\label{mutationpreserves}

Let  $\epsilon = (E_1,\dots, E_n)$ be a  complete exceptional sequences in
$\coh \XX$ and  $\sigma$ be the generator of the braid group $B_n$
such that  $\sigma.\epsilon = (E_1, \dots, E_{k-1}, L E_{k+1}, E_{k}, E_{k+2}, \dots, E_{n})$,
where we write shortly  $L E_{k+1}$ instead of $L_{E_k} E_{k+1}$.
  Then the respective dimensions
$\ddim  \Hom(E_{i}, L E_{k+1})$, $\ddim  \Ext^1(E_{i}, L E_{k+1})$ for $1 \leq i \leq k-1$,  $\ddim
\Hom_{\XX}( L E_{k+1},E_{i})$, $\ddim(  \Ext^1(L E_{k+1}, E_{i})$ for $i \in \{k, k+2, \dots, n \}$ and $rank(L
E_{k+1})$ depend only on the dimensions of the $\Hom, \Ext^{1}$ and the $ranks$ of the elements in $\epsilon$.
\end{Lem}

\Proof
In remarks \ref{dimension} and \ref{dimensionext} we have seen  that the dimensions $\ddim  \Hom( L E_{k+1},E_{k})$,
    $\ddim \Ext^1(L E_{k+1}, E_{k})$ depend only of the dimension of $\Hom_{\XX}(X_{k}, X_{k+1})$  or
$  \Ext^1(X_{k}, X_{k+1})$. Now we will prove the claim for
$\ddim  \Hom(E_{j}, L E_{k+1})$, $\ddim \Ext^1(E_{j}, L E_{k+1})$ for $1 \leq j \leq k-1$.

Suppose that the mutation is given by  type $(\alpha)$, then we have the exact sequence
$$ 0 \lrarrow L E_{k+1}  \lrarrow  \Hom(E_k, E_{k+1}) \otimes_k E_{k} \stackrel{{\mathrm {can}}}{\lrarrow}
E_{k+1} \lrarrow 0,$$
which induces a long exact sequence
\[
  \begin{array}{rcll}
    0 \lrarrow  \Hom(E_{j}  , L E_{k+1} ) &  \lrarrow   \Hom(E_{j}, \Hom(E_k, E_{k+1}) \otimes_k E_{k})  &
    \lrarrow   \Hom(E_{j}  , E_{k+1})  \lrarrow  \\
    \lrarrow   \Ext^1(E_{j}  , L E_{k+1} )  & \lrarrow
  \Ext^1(E_{j}, \Hom(E_k, E_{k+1}) \otimes_k E_{k}) & \lrarrow     \Ext^1(E_{j}  , E_{k+1})   \lrarrow 0 \\
  \end{array}
  \]
  for $1 \leq j \leq k-1$.

  Since by \cite[Lemma 3.2.4]{Meltzer:1995}
   $\Hom_{\XX}(E_{j}, E_{k+1}) = 0$ or $  \Ext^1(E_{j}, E_{k+1}) = 0$  we have either

  $$\ddim  \Hom(E_{j}, L E_{k+1}) = \ddim  \Hom(E_{j}, E_{k})  \cdot \ddim  \Hom( E_{k}, E_{k+1})
  $$
  and
  $$\ddim   \Ext^1( E_{j}, L E_{k+1}) = \ddim   \Ext^1(E_{j}, E_{k}) \cdot  \ddim  \Hom(E_{k}, E_{k+1}) -
  \ddim   \Ext^1( E_{j}, E_{k+1})$$
    or
   $$
  \ddim  \Hom(E_{j}, E_{k}) \cdot \ddim  \Hom( E_{k}, E_{k+1}) + \ddim   \Ext^1(E_{j}, L E_{k+1}) =
  $$
  $$
  \ddim  \Hom( E_{j}, L E_{k+1}) + \ddim  \Hom( E_{j}, E_{k+1}).
$$

Since $(E_{j}, L E_{k+1})$ is an exceptional pair, we have $\Hom_{\XX}( E_{j}, L E_{k+1} ) = 0$ or
$  \Ext^1(E_{j}, L E_{k+1}) = 0$.  Each one gives us that $\ddim  \Hom( E_{j}, L E_{k+1})$ and
$  \Ext^1(E_{j}, L E_{k+1})$ depend only of the dimensions of the $\Hom$, $\Ext$ spaces  of $\epsilon$.

In the cases that the left mutation is given by  type $(\beta)$ or type  $(\gamma)$ the proof is similar.
\ende
\vspace{0.2cm}
We have as a consequence of the previous discussion the following:

\begin{Cor}\label{satistythesamedimensionformulas}
Suppose that
$\epsilon = (E_1,\dots, E_n)$  and $\epsilon' = (E^{'}_1,\dots, E^{'}_n)$ are  complete exceptional sequences in
$\coh \XX$ such that the following formulas are valid
$\ddim \Hom_{\XX}(E_{j}, E_{l}) = \ddim \Hom_{\XX}(E^{'}_{j}, E^{'}_{l})$,  $\ddim\Ext^{1}_{\XX}(E_{j},  E_{l}) =
\ddim\Ext^{1}_{\XX}(E^{'}_{j},  E^{'}_{l})$
and $\rk \ (E_j) = \rk ( E^{'}_{j})$ for all $1 \leq j,l \leq n$.  Let $\sigma \in B_{n}$ and
$\sigma \epsilon = (F_{1}, \cdots, F_{n})$, $\sigma \epsilon' = (F^{'}_{1}, \cdots, F^{'}_{n})$. Then $\ddim
\Hom_{\XX}(F_{j}, F_{l}) = \ddim \Hom_{\XX}(F^{'}_{j}, F^{'}_{l})$,  $\ddim\Ext^{1}_{\XX}(F_{j},  F_{l}) =
\ddim\Ext^{1}_{\XX}(F^{'}_{j},  F^{'}_{l})$
and \\ $\rk ( F_j ) = \rk ( F^{'}_{j})$ for all $1 \leq j,l \leq n$.
\end{Cor}
 \ende

 \begin{Thm} \label{strongest}
Let $\XX=(\pp,\lala)$ and $\XX'=(\pp,\lala')$ be weighted projective lines with the same weight type.
Then  $\stgldim \XX =\stgldim \XX' $.
\end{Thm}

\Proof
Let $m$ be maximal such that there exists a tilting complex ${T}$ of the form
 $\bigoplus_{i=0}^m T_i [i]$ with $T_i \in \coh\XX$ and $T_0 \neq 0 \neq T_m$.
 Write $T = \bigoplus E_j[n_j]$ with indecomposable sheaves
  $E_j$ and $n_j \in \ZZ$.
 The $E_j$ can be ordered in such a way  that they form a full exceptional sequence
  $\epsilon $ in $\coh \XX$.
 By Theorem 1.1. there exists a braid group element $\si \in B_n$ such that $\epsilon= \si \cdot \kappa$
 where
 $\kappa=(\Oo_{\XX},  \Oo_{\XX}(\vx_{1}), \dots,  \Oo_{\XX}((p_1-1) \vx_{1}), \dots,
 \Oo_{\XX}(\vx_{t}), \dots,  \Oo_{\XX}((p_t-1) \vx_{t}), \Oo_{\XX}(\vc) )$
  is the exceptional sequence obtained from the canonical tilting sheaf
 $ \bigoplus_{0 \leq \vx \leq \vc} \Oo_{\XX}(\vx)$
  on $\XX$.

Now the application of the same braid group element $\si$ to the  exceptional sequence
 $\kappa'=(\Oo_{\XX'},  \Oo_{\XX'}(\vx_{1}), \dots,  \Oo_{\XX'}((p_1-1) \vx_{1}), \dots,
 \Oo_{\XX'}(\vx_{t}), \dots,  \Oo_{\XX'}((p_t-1) \vx_{t}), \Oo_{\XX'}(\vc) )$
obtained from the canonical tilting sheaf
 $ \bigoplus_{0 \leq \vx \leq \vc} \Oo_{\XX'}(\vx)$
  on $\XX'$ yields a full exceptional sequence $\epsilon'$ for the weighted projective line $\XX'$.

  The exceptional sheaves $\Oo_{\XX}(u \vx_{i})$
and $\Oo_{\XX'}(s \vx_{i}) $ of $\kappa$ and $\kappa'$, respectively, satisfy the same dimension
  for the $\Hom$ and $\Ext$ spaces
that is,
 $$\dim \Ext^{k}_{\XX}( \Oo_{\XX}(m \vx_{i}) , \Oo_{\XX}(n \vx_{j})) = \dim \Ext^{k}_{\XX'}(\Oo_{\XX'}(m \vx_{i}) ,
 \Oo_{\XX'}(n \vx_{j}) ) $$ for $i,j \in \{1, \cdots, t \}$, $ m,n  \in \{1, \cdots, max\{ p_{1}, \cdots, p_{t} \} \}$
and $k \in \{0, 1 \}$.
Moreover, the ranks of all sheaves of $\kappa$ and $\kappa'$ equal $1$.

The sequence $\epsilon'$ is constructed from $\kappa'$ using successively the same kind of mutations
as in the construction of $\epsilon$ from $\kappa$. Therefore, following Corolary
\ref{satistythesamedimensionformulas}
the exceptional sheaves $E'_j$ of
 $\epsilon'$ satisfy the same dimension formulas for  the $\Hom$, $\Ext$  and $rank$ spaces as
the exceptional sheaves $E_j$ of $\epsilon$.

Therefore the exceptional sheaves $E'_j$ of
 $\epsilon'$ satisfy the same dimension formulas for  the $\Hom$ and $\Ext^1$ spaces as
the exceptional sheaves $E_j$ of $\epsilon$. Hence the $E'_j$ can be shifted in the derived category of
$\coh {\XX'}$ as the $E_j$ which yields a tilting complex
 $\bigoplus_{i=0}^m T'_i [i]$ with  $T'_i \in \coh{\XX'}$ for $\XX'$ and with $T'_0\neq 0 \neq T'_m$.
 Consequently  $\stgldim \XX =m  \leq \stgldim \XX'$. By
 symmetry,  $ \stgldim \XX' \leq \stgldim \XX$
 and consequently  $ \stgldim \XX =  \stgldim \XX'$.

 \ende

\vspace{0,2cm}
Note that from Corollary \ref{satistythesamedimensionformulas} we obtain that the
ordinary quivers of the algebras  $\End \ T$ and $\End \ T'$ in the former theorem are the same
which was already stated in \cite{Meltzer:2004}.

  The former proof also suggests the following:\\
\noindent{\bf Conjecture:}
Let  ${T}$ be a  tilting complex of the form
$\bigoplus_{i=0}^m T_i [i]$ with $T_i \in \coh\XX$ and $T_0 \neq 0 \neq T_m$
and $A = \End \ T$. The strong global dimension of $A$ does not depend on the  parameters.

The validity of this conjecture implies the statement of the Theorem
\eqref{strongest}.

\section{Determinants for exceptional sequences} \label{determinants}

Let $f_1, \dots f_n$ be  group homomorphisms defined on the Grothendieck group $K_0(\XX)$ of a weighted projective line
with values in $\ZZ$. For a full exceptional sequence $\epsilon=(E_1, \dots, E_n)$ on $\XX$  we  form the $n \times n$
matrix
$M(\epsilon)$ whose coefficient at the place $(i,j)$ equals $f_i(E_j)$ and we consider the determinant of that matrix
$\det(M(\epsilon))$.

\begin{Thm}\label{det1}
There exists a constant $c \in k$ such that $\det(M(\epsilon))=c $ or $ -c$ for all full exceptional sequences
$\epsilon$
in $\coh\XX$.
\end{Thm}

\Proof
We are going to show that the determinant of the matrix does not change if we apply to the exceptional sequence
in $\coh \XX$ the left
mutation $\si_i$ . For  right mutations the proof is analogous.

For a full exceptional sequence $\epsilon = (E_1, E_2, \dots E_n)$ we denote
$\dim_k\Hom(E_i, E_{i+1})=h $ and $\dim _k\Ext^1(E_i, E_{i+1})=e$.
Now, $\si_i \cdot \epsilon $  equals  $ ( E_1, \dots, E_{i-1}, L E_{i+1}, E_i,  E_{i+2}, \dots E_n)$ and we have
$[ L_{E_i}E_{i+1}]= h[E_i]-[ E_{i+1}]$,
 $ [ L_{E_i}E_{i+1}]= [ E_{i+1}]-h[E_i]$
  or $[ L_{E_i}E_{i+1}]= e[E_i]+[ E_{i+1}]$
depending on the type of the left mutation of the pair
$(E_i, E_{i+1})$ (see section \ref{preliminaries}).
 The matrix for the exceptional sequence  $\si_i \cdot \epsilon $
is obtained from that of $\epsilon$ by replacing  the values in the $i$-th column
by $f_j( E_{i+1})-hf_j(E_i)$,
 $-f_j( E_{i+1})+hf_j(E_i)$
  or  $f_j(E_{i+1})+ef_j(E_i)$,
   $j=1,\dots, n$  and  by replacing  the values
  in the $i+1$-th column  by $f_j(E_i)$, $j=1, \dots n$.
   Then the statement follows from the known rules for determinants.
  \ende

\vspace{0,3cm}
In particular we can apply the method above to the rank function, the degree function
and the $n-2$ Euler forms $\langle -, S_{i,j} \rangle$, $j=1, \dots p_i-1$, $i=1,\dots t$.

 \begin{Cor}\label{det2}

 For each full exceptional sequence $\epsilon=(E_1, E_2, \dots, E_n) $ in $\coh \XX$
  the determinant of the matrix
\begin{footnotesize}
\begin{equation*}
M(\epsilon)=
\begin{pmatrix}
\rk E_1 & \rk E_2 & \rk E_3 & \dots & \rk E_{p_1+1}& \dots & \rk E_{n} \\
\deg E_1 & \deg E_2 & \deg E_3 & \dots & \deg E_{p_1+1} &\dots  &  \deg E_{n} \\
\langle E_1,S_{1,1}\rangle& \langle E_2,S_{1,1}\rangle & \langle E_3,S_{1,1} \rangle &\dots &\langle
E_{p_1+1},S_{1,1}\rangle & \dots  &\langle E_n,S_{p_1,1}\rangle\\
\vdots & & & & & & &  \\
\langle E_1,S_{1,p_1-1}\rangle& \langle E_2,S_{1,p_1-1}\rangle & \langle E_3,S_{1,p_1-1} \rangle &\dots &\langle
E_{p_1+1},S_{1,p_1-1}\rangle & \dots &\langle E_n,S_{1,p_1-1}\rangle\\
\vdots & & & & & & &  \\
\langle E_1,S_{t,1}\rangle& \langle E_2,S_{t,1}\rangle & \langle E_3,S_{t,1} \rangle &\dots &\langle
E_{p_1+1},S_{t,1}\rangle & \dots  &\langle E_n,S_{t,1}\rangle\\
\vdots & & & & & & &  \\
\langle E_1,S_{t,p_t-1}\rangle& \langle E_2,S_{t,p_t-1}\rangle & \langle E_3,S_{t,p_t-1} \rangle &\dots &\langle
E_{p_1+1},S_{t,p_t-1}\rangle & \dots &\langle E_n,S_{t,p_t-1}\rangle\\
\end{pmatrix}
\end{equation*}
\end{footnotesize}
equals $p$ or $-p$. Recall that $p$ denotes the least common multiple of the weights $p_1,\dots,p_t$.
\end{Cor}

\Proof
The determinant is easily  calculated to be $p$ or $-p$ for the exceptional sequence \\
$(\Oo, \Oo(\vc), S_{1,1}, \dots, S_{1,p_1-1}, \dots S_{t,1}, \dots, S_{t,p_1-1})$
using the block structure of the matrix and the fact that $\rk \Oo= \rk \Oo(\vc)=1$,
$\deg \Oo= 0$ and $\deg \Oo(\vc)= p$.
Then the statement follows from Theorem \ref{det1}.
\ende

\vspace{0,3cm}
\begin{Rem}
The determinantal equation obtained in the way above can be interpreted as a diophantine equation
for the weighted projective line $\XX$.
Diophantine  equations expressed for data in terms of exceptional sequences seems to be typical.
So  Rudakov showed that the ranks of the vector bundles of an exceptional triple on the projective
plane satisfy the Markov equation $X^2+Y^2+Z^1=3XYZ$ \cite{Ru}.
Diophantine equations for partial tilting sequences on weighted projective lines were given in
\cite[Chapter 10.2]{Meltzer:2004}.

\end{Rem}

\vspace{ 2cm}

\begin{center}{Addresses}
\end{center}


\begin{thebibliography}{}

\bibitem[ALM]{ALM}
  E. R. Alvares, P. Le Meur, and E. N. Marcos. \emph{The strong global dimension of piecewise hereditary algebras}. J. Algebra, 481:36-67, 2017.

\bibitem[S]{S}
S. C. Schmidt,
\emph{Complexos Tilting e Dimens\~ao Global Forte em Álgebras Heredit\'arias por partes.} PhD Thesis. Universidade Federal do Paran\'a 2017.



\bibitem[B]{B}
A. I. Bondal,
 \emph{Representation of associative algebras and coherent sheaves},
 Math. USSR, Izv. 34, No. 1, 23-42 (1990); translation from Izv. Akad. Nauk SSSR, Ser. Mat. 53, No. 1, 25-44 (1989).

\bibitem[CB]{CB}
W. Crawley-Boevey,
 \emph{Exceptional sequences of representations of quivers},
 Dlab, Vlastimil (ed.) et al., Representations of algebras. Proceedings of the sixth international conference on
 representations of algebras, Carleton University, Ottawa, Ontario, Canada, August 19-22, 1992. Providence, RI:
 American Mathematical Society. CMS Conf. Proc. 14, 117-124 (1993).

\bibitem[GL1]{Geigle:Lenzing:1987}
W. Geigle and H. Lenzing, \emph{A class of weighted projective curves
  arising in representation theory of finite-dimensional algebras},
  Singularities, representation of algebras, and vector bundles ({L}ambrecht,
  1985), Lecture Notes in Math., vol. 1273, 265-297 Springer, Berlin, (1987).



\bibitem[GL2]{Geigle:Lenzing:1991}
W. Geigle and H. Lenzing, \emph{Perpendicular categories with applications
to representations and sheaves},
  J. Algebra 144, No. 2, 273-343 (1991).

\bibitem[GR]{GR}
A. L. Gorodentsev and A. N. Rudakov,
\emph{Exceptional vector bundles on projective spaces},
Duke Math. J. 54, 115-130 (1987).

\bibitem[H]{Ha} D. Happel,
\emph{Perpendicular categories to exceptional modules},
An. Stiint. Univ.  “Ovidius” Constanta, Ser. Mat. 4, No. 2, 66-75 (1996).

\bibitem[HR]{HR} D. Happel and  C. M. Ringel
\emph{Tilted algebras},
Trans. Am. Math. Soc. 274, 399-443 (1982).





\bibitem[LP]{LP}
H. Lenzing and J. A. de la Pe\~{n}a, \emph{Wild canonical algebras},
Math. Z. 224, No. 3, 403-425 (1997).


\bibitem[M1]{Meltzer:1995}
H. Meltzer, \emph{Exceptional sequences for canonical algebras},
Arch. Math. 64, No. 4, 304-312 (1995).

\bibitem[M2]{Meltzer:2004}
H. Meltzer, \emph{Exceptional vector bundles, tilting sheaves and tilting
  complexes for weighted projective lines},
  Mem. Am. Math. Soc. 808, 138 p. (2004).

\bibitem[Ri]{Ri}
C. M. Ringel, \emph{Tame algebras and integral quadratic forms},
Lecture Notes in Mathematics. 1099. Berlin etc.: Springer-Verlag. XIII, 376 p. (1984).


\bibitem[Ru]{Ru} A. N. Rudakov,
 \emph{The Markov numbers and exceptional bundles on $\PP^2$},
Math. USSR, Izv. 32, No. 1, 99-112 (1989); translation from Izv. Akad. Nauk SSSR, Ser. Mat. 52, No. 1, 100-112 (1988).

\end{thebibliography}
\end{document}